\documentclass[runningheads]{llncs}

	\usepackage[T1]{fontenc}
	\usepackage{graphicx}
	\usepackage{amssymb}
	\usepackage{amsmath}

	\usepackage[numbers]{natbib}
	
	\usepackage{float}
	\usepackage{multirow}
	\usepackage{xcolor}
	\usepackage{forest}
	\usepackage[ruled]{algorithm2e}

	\usepackage[refeq]{myMac-els}

\newcommand{\icmsOrArxiv}[2]{#1}	

\newcommand{\whh}{{\wh{h}}}

\newcommand{\End}{\mbox{End}} 
\newcommand{\Image}{\mbox{Image}} %
\newcommand{\EndCover}{\mbox{EndCover}} 

\newcommand{\Tube}{\mbox{Tube}} 
\newcommand{\NoTube}{{\bf NoT}} 
\newcommand{\withT}{{\bf T}} 
\newcommand{\NoB}{{\bf NoB}} 
\newcommand{\withB}{{\bf B}} 
\newcommand{\wmax}{w_{\max}} 
 
\newcommand{\mumax}{\mu_{\max}} 

\newcommand{\nnorm}[1][\cdot]%
	{\left\lVert #1 \right\rVert_{\max}}

\def\lognorm{{logNorm}}

\def\ivp{\mbox{{\ttt{IVP}}}}		
\def\stepA{{\mbox{\tt StepA}}}
\def\stepB{{\mbox{\tt StepB}}}

\def\Direct{\mbox{\tt Direct}}
\def\PointError{\mbox{\tt PE}}
\def\RangeError{\mbox{\tt RE}}

\def\Refine{\mbox{\ttt{Refine}}}
\def\Extend{\mbox{\ttt{Extend}}}
\def\Split{\mbox{\ttt{Split}}} 

\def\taylortube{{\ttt{TaylorTube}}}
\def\Bisect{{\ttt{Bisect}}}

\def\stage{{\calS}} 

	\newcommand{\taylor}{^{\text{taylor}}}
	\newcommand{\Tay}[2][p]{ T^{#1}_{#2}} 

\renewcommand{\icmsOrArxiv}[2]{#1}	

\AtBeginDocument{%
  }

	\def\tsc#1{\csdef{#1}{\textsc{\lowercase{#1}}\xspace}}
	\tsc{WGM}
	\tsc{QE}
	\tsc{EP}
	\tsc{PMS}
	\tsc{BEC}
	\tsc{DE}
\newcommand{\pmcell}[2]{%
  \shortstack[c]{$#1$\\[1pt]$\pm #2$}%
}
\begin{document}
\icmsOrArxiv{
}{ 
}

\title{Taylor Tube Method for Validated IVP\thanks{
	This document contains the results of the research funded by NSF
	Grant \#CCF-2212462.
    \\
    A paper by the same title will appear in the Proceeding of the International Congress on Mathematical Software (ICMS) 2026. This arxiv version has all the proofs in an appendix.}}

\titlerunning{Taylor Tube}

\author{
Bingwei Zhang\inst{1}\orcidID{0009-0002-2619-9807}
\and
Chee Yap\inst{1}\orcidID{0000-0003-2952-3545}
}
\authorrunning{B. Zhang and C. Yap}

\institute{
Department of Computer Science\\
Courant Institute School of Mathematics, Computing, and Data Science\\
New York University, New York, USA\\
\email{bz2517@nyu.edu, yap@cs.nyu.edu}
}
\maketitle
\begin{abstract}
	We recently introduced a novel architecture for
	the design of validated IVP algorithms.
	This architecture forms the basis of our complete
	validated algorithm for IVP.
	A key subroutine in our algorithm is
	the \dt{Euler Tube}: it gave a technique for refining
	end- and full-enclosures and is also key to deriving
	a complexity bound of our IVP solver.
	In this paper, we generalize it to
	\dt{Taylor Tube} of degree $p\ge 1$.
	As expected, higher-degree
	Taylor Tubes improve accuracy. But surprisingly,
	our experiments show that it can
	also lead to an overall speedup
	when combined with bisection.

\end{abstract}


\sect{Introduction}

Consider the initial value problem (IVP)
for an autonomous ODE of the form
	\beql{bfx'}
		\bfx'=\bff(\bfx).\eeql
where $\bff:\RR^n\to\RR^n$ is $C^k$-continuous (for some
global constant $k\ge 1$).
For $B_0\ib\RR^n$ ($n\ge 1$) and $h>0$, let
	$\ivp_\bff(B_0,h)$ denote the set of solutions
$\bfx:[0,h]\to \RR^n$ that satisfies \refeq{bfx'} and
$\bfx(0)\in B_0$. 
We say $(B_0,h)$ is \dt{valid} (or $\bff$-valid)
if for each $\bfp_0\in B_0$,
some $\bfx\in \ivp_\bff(B_0,h)$ with $\bfx(0)=\bfp_0$
has unique-existence.\footnote{
	i.e., it exists, and it is unique. }

For a valid $(B_0,h)$, let
	\beql{endImage}
		\mmatX[rcl]{
		\End_\bff(B_0,h)
			&\as& \set{\bfx(h): \bfx\in \ivp_\bff(B_0,h)},\\
		\Image_\bff(B_0,h)
			&\as& \set{\bfx(t): \bfx\in \ivp_\bff(B_0,h), t\in [0,h]}.
		}\eeql
We fix $\bff$ in this paper; so the subscript $\bff$
may be omitted above.
By an \dt{end-closure} and a \dt{full-enclosure}
for $(B_0,h)$ we mean any sets $E$ and $F$ that contain
$\End(B_0,h)$ and $\Image(B_0,h)$, respectively.
Two central problems in the applications of IVP are
to compute such end- and full-enclosures.
These are often called \dt{reachability problems} in control theory.
The literature on validated algorithms for IVP has a history of
over 50 years going back to Moore, and are of course
based on interval methods.
See the surveys of Corliss \cite{corliss:survey-ode-intvl:89}
and Nedialkov et al.~\cite{nedialkov+2:validated-ode:99}).

\ignore{%
In \dt{validated algorithms}, any 
computed numerical quantity or geometric set must be
explicitly contained in an output interval, box or ball.
}

We formulate the following reachability problem, called the
\dt{End Cover Problem}:
Given $(\bff, B_0, h, \veps)$ ($\veps>0$) where $B_0$ is a box,
to compute a finite set of boxes $\calC=\set{B_1\dd B_m}$
such that
	\beql{endbff}
		\End_\bff(B_0,h)~ \ib~ \bigcup\calC
			~\ib~ \End_\bff(B_0,h) \oplus [-\veps,\veps]^n\eeql
where
$\bigcup\calC=\cup_{i=1}^m B_i$ and
$A \oplus B$ denotes the Minkowski sum of Euclidean sets $A,B\ib\RR^n$.
Any $\calC$ that satisfies \refeq{endbff} is called an
\dt{$\veps$-cover} of $\End_\bff(B_0,h)$.

	\renewcommand{\alt}[2]{#2} 
	\renewcommand{\alt}[2]{} 

	\alt{
		\beql{endCover}
			\includegraphics[width=0.8\columnwidth]{figs/endCoverProb}
		\eeql
	}{
		\Ldent\progb{
			\lline[0]
					\mbox{\sc End Cover Problem:}
			\lline[10]
					$\EndCover_\bff(B_0,\veps; \coblue{h=1}) \to \calC$
			\lline \myhlineY
			\lline[0] INPUT: \quad $B_0\in\intbox\RR^n$,
						 $\veps>0$, \coblue{$h>0$}
				where $(B_0,h)$ is valid.
			\lline[0] OUTPUT: $\calC$ is
					an \dt{$\veps$-cover} of $\End_\bff(B_0,h)$.
			\lline[15] 
					I.e., $\calC$ is a finite set of boxes such that
			\lline[10] 
					$\End_\bff(B_0,h) ~\ib~ \bigcup\calC ~\ib~ 
						\Big(\End_\bff(B_0,h) \oplus [\pm \veps]^n\Big) $
		}
	}

Recently, we introduced an algorithm
\cite{zhang-yap:ivp:25arxiv,zhang-yap:cover:26arxiv}	
for this End Cover Problem.  It represented the
first \dt{complete} validated IVP algorithm.
By ``complete'' we mean that our algorithm is a halting algorithm.
Most validated algorithms are only \dt{partially correct}\footnote{
	This is a standard terminology in theoretical computer science.} 
in the sense that, the output is
correct {\em if it halts}.  Halting is especially challenging
in IVP problems, partly because the error typically
grows exponentially with time.
We note some features in our solution:
(1) We guarantee that the ``over-approximation'' of $\calC$
is within any user-specified $\veps$ bound.  To our knowledge,
no current solver promises such bounds.
(2) Our algorithm is fully automatic; it does not require \dt{hyper-parameters}, such as a minimum step size or tolerance bounds, which are commonly used in such algorithms.
(3) In preliminary implementations, our algorithm 
has proved to be practical. 
It is reasonable to expect that the algorithm takes a performance hit
in order to ensure completeness.  But surprisingly, it
often performs better than the incomplete algorithms.

Our previous papers  \cite{zhang-yap:ivp:25arxiv,zhang-yap:cover:26arxiv}
introduced a novel architecture for
IVP algorithms based on the \dt{scaffold} framework.
Within this framework, many of our original
subroutines could be swapped.  One of these 
subroutines is called \dt{Euler tube}.  In this paper, we
generalize it into \dt{Taylor tube}, based on
using the first $p\ge 1$ terms of a Taylor expansion.
The case $p=1$ corresponds to Euler tube.
Our experiments show that a higher degree for $p$
leads to tighter enclosures (not a surprise) but
surprisingly, they can also lead to overall speedup.
%
We remark that such tube techniques were critical
for our ability to prove complexity bounds for our
halting algorithm.

\section{Taylor Tube}
	Most of the following notations are taken
	from \cite{zhang-yap:ivp:25arxiv,zhang-yap:cover:26arxiv}.
	For a square matrix $A\in\RR^{n\times n}$ and the induced matrix
	$p$-norm, the \dt{logarithmic $p$-norm} is
		\[
		\mu_p(A)\as \lim_{h\to 0+}\frac{\|I+hA\|_p-1}{h}.
		\]
	We focus on $p=2$ and write $\mu(A)$ as the \lognorm.
	For any $B\ib\RR^n$, define
		\[
		\mu(J_\bff(B))\as \sup\set{\mu(J_\bff(\bfp)):\bfp\in B}.
		\]
	Call $\olmu$ a \dt{\lognorm-bound} on $B$
	if $\olmu \ge \mu(J_\bff(B))$.
	Unlike classical norms, logarithmic norms can be negative.
	We have this basic result:

\ignore{%
	The following is from \cite{neumaier:rigorous-bds:93}:
	
	\bthmT[ne]{Fundamental IVP Inequality}
	Let $\bfx\in \ivp_\bff(\bfp_0;h)$ and
	$\xi\in C^1([0,h]\to\RR^n)$ be any approximate trajectory.
	Assume constants
	$\veps,\delta,\olmu\ge 0$ satisfy
	\begin{enumerate}
		\item
		$\displaystyle
		\veps \ge
		\|\xi'(t)-\bff(\xi(t))\|
		\quad \text{for} all\ t\in[0,h],$
		\item
		$\displaystyle
		\delta \ge
		\|\xi(0)-\bfx(0)\|,$
		\item
		$\displaystyle
		\olmu \ge
		\mu\!\big(
		J_\bff(s\bfx(t)+(1-s)\xi(t))
		\big)$
		for all $s\in[0,1]$ and $t\in[0,h]$.
	\end{enumerate}
	Then for all $t\in[0,h]$,
	\beql{xibfx}
	\|(\xi(t)-\bfx(t))\|
	\le
	\begin{cases}
		\delta e^{\olmu t}
		+\dfrac{\veps}{\olmu}(e^{\olmu t}-1),
		& \olmu\ne 0,\\[1ex]
		\delta+\veps t,
		& \olmu=0.
	\end{cases}
	\eeql
	\ethmT
}%
	\bleml[FundIneq]
	Let $\bfx_1,\bfx_2\in \ivp(Ball_{\bfp_0}(r),h,F)$ and assume
	$\olmu\ge \mu_2(J_\bff(F))$.  For all $t\in[0,h]$,
		\beql{trajSep}
			\|\bfx_1(t)-\bfx_2(t)\|_2 \le
			\|\bfx_1(0)-\bfx_2(0)\| e^{\olmu t}.
		\eeql
	\eleml
	
	This lemma is just a special case of
	\lognorm\ version of the \dt{Fundamental IVP Inequality}
	(\cite{neumaier:theoryI:94,hairer+2:ode1:bk}).

\medskip\noindent
	\dt{Tubular region of a trajectory.}
	Fix a trajectory
			$\bfx= \ivp(\bfq_0,h)$.
	For any $\delta>0$, the \dt{$\delta$-tube} of $\bfx$ is the set
	\[
		\Tube_\delta(\bfx)\as
		\set{(t,\bfp): 0\le t\le h,\ \|\bfp-\bfx(t)\|_2\le \delta}
		\quad (\ib [0,h]\times\RR^n).
	\]
	An arbitrary curve $C:[0,\infty) \to\RR^n$
	\dt{lies in the $\delta$-tube of $\bfx$} if
	$(t,C(t))\in \Tube_\delta(\bfx)$ for all $t\in[0,h]$.
	For $p\ge 1$, the \dt{Taylor polynomial} of $\bff$
	of degree $p$ at any $\bfq\in\RR^n$ is defined to be
		\beql{taylorpoly}
			\Tay[p]{\bfq}(t)
				\as \sum_{i=0}^{p} t^i \bff\supn[i](\bfq))
		\eeql
	where 
	\footnote{
		We may call this the \dt{normalized
		$i$-th Taylor differential coefficient} of $\bff$,
		where ``differential'' here emphasizes
		that it is not the usual coefficient of polynomial.
	}
		$$\bff\supn[i](\bfx)\as  \clauses{
				\bfx &\rmif\ i=0,\\
						\tfrac{1}{i} \left(
				J_{\bff\subn[i-1]}\Bigcdot \bff\right)(\bfx)
						&\elsE}
		$$
	where $J_\bfg$ denotes the Jacobian of $\bfg$
	and $\Bigcdot$ is matrix-vector product.
	The \dt{Taylor sequence} of $\bfq_0$ with \dt{step size} $\olh >0$ 
		is the infinite sequence of points
			$\bfq_0, \bfq_1, \bfq_2 \ldots$
		where 
			$\bfq_{i+1} \as \Tay[p]{\bfq_i}(\olh)$ for all
			$i\ge 0$.  
	This sequence defines the piece-wise
	continuous \dt{Taylor curve} where
		\beql{taylorcurve}
			\Tay[p]{\bfq_0,\olh}(t) \as
					\Tay[p]{\bfq_i}(t-i\olh),
					\qquad (t\ge 0, i=\floor{t/\olh}).
			\eeql

	
	We could similarly define the
	\dt{$\delta$-tube of the Taylor curve}
	$\Tay[p]{\bfq_0,\olh}(t)$ and
	ask: does it contain the trajectory $\bfx$?
	Our next lemma gives an affirmative answer
	provided the step size $\olh$ is less than an explicit bound.
	Using \refLem{FundIneq}, we obtain:
	
	\blemT[taylorStep]{Taylor Tube Step Size Bound}
	Let $\bfx =\ivp(\bfq_0,h)$ and $B_1$ be a full enclosure
	of $(\bfq_0,h)$. 
	Assume that $\olmu$ is a \lognorm-bound on $B_1$
		and $\olM \ge \|\bff\supn[p+1](B_1)\|$.
	For any $\delta>0$, if $\olh$ satisfies 
		{\small
	    \beql{h1}
			\olh  \le h\taylor_p(h,\olM,\ol\mu,\delta)
			\as \begin{cases}
		\displaystyle  \Big(
			\dfrac{\olmu\,\delta }{\olM\,(e^{\ol\mu h}-1)}
			\Big)^{\!1/p} & \text{if }\ol\mu> 0,\\[10pt]
		\displaystyle  \Big(
			\dfrac{\delta }{\olM h}
			\Big)^{\!1/p} & \text{if }\ol\mu= 0,\\[10pt]
		\displaystyle 
			\min\!\left\{ \left(\frac{\olmu\delta}{2\olM\,(e^{\olmu h}-1)}
\right)^{1/p},\, -\frac{1}{\olmu} \right\}, & \text{if }\ol\mu<0,
		\end{cases}
	\eeql}
and if $\Image(\bfq_i,\olh)\ib B_1$ (for all $i=0\dd \floor{h/\olh}$ ) then
        the $\delta$-tube of $\bfx$ contains the Taylor curve 
			$\Tay[p]{\bfq_0,\olh}(t)$ of step size $\olh$.
	\elemT

	\ignore{
	Note that we could also write $\Tay[p]{\bfx,\olh}(t)$
	as $\Tay[p]{\bfq_0,\olh}(t)$ since $\bfx$ is the unique
	trajectory with $\bfx(0)=\bfq_0$.
	}%

\sect{Overview of the End-Cover Algorithm}
	Most validated IVP algorithms are based on two key subroutines
	which we call\footnote{
		The $\stepA$ in \cite{zhang-yap:ivp:25arxiv},
		has two extra arguments $h,\veps\ge 0$, which is
		defaulted to $h=\infty, \veps=0$ here.
	}
	\stepA\ and \stepB.  Given a box $E_0$,
	we may call \stepA/\stepB\ in sequence to transform $E_0$ into
	the following triple and then quad:
		\beql{seqAB}
			E_0\longder[\stepA] (E_0,h,F_1)
				\longder[\stepB] (E_0,h,F_1,E_1).
			\eeql
	In general, for $i\ge 1$,
	we transform $E_{i-1}\der (E_{i-1},h_i,F_i,E_i)$.
	These triples and quads are said to be \dt{admissible}
	and have the property that $F_i$ and $E_i$ are
	full- and end-enclosures for $(E_{i-1},h_i)$.
	Note that $E_m$ ($m\ge 1$) is an end-enclosure
	for $(E_0, t_m)$ where $t_m=\sum_{i=1}^m h_i$.
	
	Next, we introduce the \dt{scaffold data structure}
	$\calS$ to store such quads.  If there
	are $m$ quads, we use these dot-notations
	for the various parameters in these quads:
		$$\calS.m,\quad \calS.E_0, \quad\calS.t, \quad\calS.E_m$$
	for the number of quads, the initial box of
	the first quad, the total time $t_m$, and end-enclosure
	at time $t_m$.  
	Each quad represents a \dt{stage} of $\calS$.
	The goal is to ensure that
		(a) $\calS.t =H$ and
		(b) $width(\calS.E_m)<\veps$
	where $H,\veps$ are the input arguments.
	If (a) fails, we call a subroutine $\calS.\Extend(H)$
	that adds a new stage (using \stepA/\stepB).
	If (b) fails, we call $\calS.\Refine(\veps)$ to 
    reduce $width(\calS.E_m)$.
	However, condition (b) may not be possible without subdividing
	$\calS.E_0$ into subboxes.  After subdivision, each subbox $E'_0\ib E_0$
	becomes the initial box of a new $m$-stage scaffold
	that inherits all the quads of $\calS$.
	The basic operation of $\calS.\Refine(\veps)$
	proceeds in \dt{phases}:
		each phase refines all $m$ stages
		(from stage $i=1\dd m$ in this order).
		Refining stage $i$ means to convert the $i$th quad
		\beql{refineStage}
			(E_{i-1},h_i,F_i,E_i)
			\der[refine] (E_{i-1},h_i,F'_i,E'_i)	\eeql
	to a new admissible quad with $F'_i\ib F_i$ and $E'_i\ib E_i$.
	If $Q_i=(E_{i-1},h_i,F_i,E_i)$ then we refine $Q_i$
	by generating $2^{\ell_i}-1$ uniform \dt{mini-steps} between
	$Q_{i-1}$ and $Q_i$:
		{\small \beql
		$Q_{i-1} \der Q_i\supx[1]
				\der Q_i\supx[2]
				\der \cdots
				\der Q_i\supx[{2^{\ell_i}-1}]
				\der Q_i.
		\eeql}
	We provide two subroutines $\Bisect$ and the
	``taylor tube method''
	for doing this refinement: $\Bisect$ simply
	calls \stepB\ twice with half the current mini-step size.
	When the mini-step size is smaller than the bound in
	\refLem{taylorStep}, we can speed up by using the
	method described in the next section.

\icmsOrArxiv{
}{ 
	See Appendix XXX for the pseudo-code of End-Enclosure Algorithm.
}

\sect{Refining End- and Full-Enclosures with Taylor Tubes}
	We now show how the Taylor Tube method can be combined
	with \stepB.  
	One of the best previous methods for \stepB\ is the Direct Method
	in \cite{nedialkov:thesis:99,nedialkov+2:validated-ode:99}:
    {\scriptsize
		\beqarrayl{direct}
		\Direct(E_0,h,F_1)
			&=& \underbracket{
					\sum_{i=0}^{k-1} h^i\bff^{[i]}(\bfm(E_0))
						}_{\text{Point}}
		+ \underbracket{
			\vphantom{\sum_{i=0}^{k-1}}h^k\bff^{[k]}(F_1)
					}_{\text{Point Error}}
		+ \underbracket{
			\Big(\sum_{i=0}^{k-1}h^iJ_{\bff^{[i]}}(E_0)\Big)
					\Bigcdot\big(E_0-\bfm(E_0)\big)
							}_{\text{Range Enclosure}}\nonumber\\
			&=& \Tay[k-1]{\bfm(E_0)}(h)
				+ \PointError(h,F_1)
				+ \RangeError(h,E_0)
		\eeqarrayl
        }
    Note that
    $ \Direct(E_0,h,F_1)$ and  $\Direct(E_0,[0,h],F_1) $
	produce (respectively) an end-enclosure ($E_1$) and
	a refinement of the full-enclosure ($F_1$) for $(E_0,h)$.
	Our next theorem provides a new method to refine $F_1$:
	
	\bthmT[taylorTube]{Taylor Tube Enclosures}\ \\
	Consider an admissible triple $(E_0, h, F_1)$ 
	where $E_0 = \mathrm{Ball}_{\bfq_0}(r_0)$.\\
	Let $\delta>0$,
		$\olh \le \min\{h, h\taylor_p(h, \olM, \olmu, \delta)\}$,
		   $\olmu$ is a \lognorm-bound on $F_1$ and $\olM\ge \|\bff\supn[p+1](F_1)\|$.
\\
	If
		\beql{image}
    \image(\bfq_i,\olh)\ib F_1 \quad(\text{for all } i=0\dd \floor{h/\olh})
    \eeql
	then the following holds:
	\benum[(a)]
	\item
		An end-enclosure for $\ivp(E_0, h)$ is given by  
			$ \mathrm{Ball}_{\Tay{\bfq_0,\olh}(h)}
				\left( r_0 e^{\olmu h} + \delta\right).  $
	\item  
		A full-enclosure for $\ivp(E_0, h)$ is given by  
		$  \Tay{\bfq_0,\olh}([0,h]) \pm [-r,r]^n,$
		where  $r = \max_{h' \in [0,h]} 
						\left(r_0 e^{\olmu h'} + \delta\right)=\delta+\max
						\{r_0, r_0e^{\olmu h} \}.$
	\ignore{
    \item  
		 If $0 \notin (\bff(F_1))_i$ for\footnote{
		 	where $( S )_i$ is the set of $i$th components
			of points in $S\ib\RR^n$.  }
		all $i=1\dd n$, then 
	    $\mathrm{Box} \!\left( E_0\cup
					\mathrm{Ball}_{\Tay{\bfq_0,\olh}(\olh)}
					\left(r_0 e^{\olmu \olh} + \delta\right) \right) $
	    is a full-enclosure for 
		$\ivp(E_0, \olh)$.
        }
	\eenum
	\ethmT
	We can combine \Direct\ with \refThm{taylorTube} to define a tighter
	end-enclosure $E_1(h)$ defined as follows:
	{\scriptsize
	\beql{E1h1}
		E_1(h) \as
		\begin{cases}
			\Direct(E_0,h,F_1)\cap
				\mathrm{Ball}_{\Tay{\bfq_0}(h)}
				\bigl(r_0 e^{\olmu h}+\delta\bigr),
				& \text{if } p\ne k-1,\\[6pt]
		\Tay[k-1]{\bfm(E_0)}(h)
			+\bigl(\PointError(h,F_1)\cap Box(\delta)\bigr)
				+\bigl(\RangeError(h,E_0)\cap Box(r_0 e^{\olmu h})\bigr),
				& \text{if } p=k-1.
		\end{cases}
	\eeql}
	Furthermore,
		$ E_1([0,h]) \cap F_1 $
	gives a tighter full-enclosure than $F_1$.

    Our Taylor Tube algorithm can 
    check \refeq{image}
    by using this inclusion
    \beql{image1}
    T_{\bfq_i,\olh}([0,\olh])+[0,\olh]^{p+1}\bff\supn[p+1](F_1)\ib F_1.
    \eeql
    Note that $ T_{\bfq_i,\olh}([0,\olh])$
    and $\bff\supn[p+1](F_1)$
    have already been computed
    when computing $E_1([0,\olh])$ \refeq{E1h1} and
    $\olM$.

\ignore{
\sect{EndCover Algorithm}
	\dt{SCAFFOLD NOTATIONS:}
		Recall that an $m$-stage scaffold
			$\stage=(\bft,\bfE,\bfF,\bfG)$ 
		has components $\bft=(t_0<t_1<\cdots<t_m)$,
			$\bfE=(E_0, E_1\dd E_m)$,
			$\bfF=(F_1, F_2\dd F_m)$,
			$\bfG=(G_1, G_2\dd G_m)$
			where 
		\beql{Gi} 
        G_i = G_i(\stage) =(
				\cored{ (\ell_i, \ol\bfE_i, \ol\bfF_i)},
				\cored{( \ol\bfmu_i, \delta_i, h\taylor_i)}).
                \eeql 
	Here we focus on the components in red.  We call
	$\cored{(\ell_i,\olE_i,\olF_i)}$ the
		$i$th \dt{mini-scaffold} and
	$\ell_i\ge 0$ is the \dt{level}.  It represents
	a subdivision of $[t_{i}-t_{i-1}]$
	into $2^{\ell_i}$ mini-steps of size
	$\whh_i\as (t_i-t_{i-1})2^{-\ell_i}$, and
	$\olE_i, \olF_i$ are $2^{\ell_i}$-vectors representing admissible
	quads $(\olE_i[j-1]), \whh_i, \olF_i[j],\olE_i[j])$
	(for $j=1\dd 2^{\ell_i}$).

	\bitem
	\item
		$\stage.m$ returns the number of stages in the scaffold.
	\item
		$\stage.t=t_m$ is the last time in $\stage$.
	\item
		$\stage.G(i).\whh$ is the ministep in $\stage_i$
		computed from data in $G(i)$.
			Give the formula for $\whh$ here...
	\item
		$\stage[i].h\taylor$ is the Taylor target step in $\stage_i$.
			Give the formula for $h\taylor$ here...
	\item
		$\stage.r$ and $\stage.R$
		are (respectively) the max width of the $E_0$ (initial)
		and $E_m$ (last) stage,
			$\stage.r\as \wmax(E_0(\stage))$,
			$\stage.R\as \wmax(E_m(\stage))$.
		$\stage.r$ is max width of the last ($m$th) stage,
			$\wmax(E_m(\stage))$.
	\item
		$\stage.\mumax$ is maximum of
			$\stage.\bfG[i].\ol\bfmu[j]$ for all
				$j=1\dd 2^{\ell_i}$ and
				all $i=1\dd \stage.m$.
	\item
		If $B\in\intbox\RR^n$, let $\Split(B)$ denote
		the set of all $2^n$ congruent subboxes of $B$
		obtained by splitting each of the $n$ dimensions of $B$.
	\item
		$\stage.\Split()\ssa Q'$ where $Q'$ is a set
		of $2^n$ scaffolds; each $\stage'$ in $Q'$
		is identical to $\stage$ except that $E_0(\stage')$
		is one of the boxes in $\Split(E_0(\stage))$.
	\eitem

{\fontsize{8}{8.6}\selectfont 
	\Ldent\progb{
		\lline[-3]
		$\EndCover_\bff(B_0,\veps_0,H)\to \calC$
		\lline[0] INPUT: \hspace*{3mm}	 $(B_0,H,\veps_0)$
		\lline[10] where $(B_0,H)$ is valid,
						$B_0\in\intbox \RR^n$ and $\veps_0,H>0$.
		\lline[0] OUTPUT: \hspace*{1mm} $\calC\ib\intbox\RR^n$
		\lline[10] that is a $\veps_0$-cover of $\End_\bff(B_0,H)$
		\lline \myhlineY
		\lline[5] Initialize a queue of scaffolds: $Q \ass\{\stage_0\}$
		\lline[10] where $\stage_0$ is a $0$-stage scaffold with
						$E_0(\stage)=B_0$
		\lline[5] Initialize a queue of boxes: $\calC \ass \es$
		\lline[5] While ($Q \neq\emptyset$)
		\lline[10] $\stage\ass Q .pop()$
					\Commentx{Get a scaffold $\stage$}
		\lline[10] While ($\stage.t<H$)
					\Commentx{$\stage.t$ is end time of $\stage$}
		\lline[15] $t\ass \stage.\Extend(H)$.
					\Commentx{$\stage.m$ has been incremented (new stage
					added)}
		\lline[15] While ($\stage.R >\veps_0$)
			\Commentx{Refine the current scaffold}
		\lline[20] For ($i=1\dd \stage.m$)
		\lline[25] If ($\stage.\bfG[i].\whh>\stage[i].h\taylor$)
			\Commentx{If ministep size $>$ Taylor target}
		\lline[30] $\stage.\Bisect(i)$
		\lline[25] Else
		\lline[30] $\stage.TaylorTube(i)$ 
					\Commentx{Call Taylor-tube on stage $i$}
		\lline[20] If $\Big(
			\stage.r\cdot e^{\stage.\mumax\cdot\stage.t} >\veps_0/2
			\Big)$
		\lline[25] $Q'\ass \stage.\Split()$
			\Commentx{$Q'$ is a queue of $2^n$ scaffolds}
		\lline[25] $\stage\ass Q'.pop()$
			\Commentx{Update $\stage$}
		\lline[25] $Q.add(Q')$
			\Commentx{This adds $2^n-1$ scaffolds to $Q$}
		\lline[15] \hspace*{0mm}End While
				\Commentx{Now $\stage.R \le\veps_0$}
		\lline[10] \hspace*{0mm}End While
				\Commentx{Now $\stage.t =H$}
		\lline[10] $\calC.append(E_{back}(\stage))$
		\lline[5] \hspace*{-3.0mm}End While
		\lline[5] Return $\calC$
	}
}

}%

\sect{Experiments}
	We present experimental evidence for the
	effectiveness of our Taylor tube techniques, using
	the list of ODE examples in \refTab{problems}.
	
	\newcommand*{\myalign}[2]{\multicolumn{1}{#1}{#2}}
	\begin{table}[] \centering
		\caption{List of IVP Problems}
		\label{tab:problems}
		{\tiny
			\btable[l| l | l | l | l | l ]{
				Eg* & \myalign{c|}{\dt{Name}}
				& \myalign{c|}{$\bff(\bfx)$}
				& \myalign{c|}{\dt{Parameters}}
				& \myalign{c|}{\dt{Box} $B_0$}
				& \myalign{c}{\dt{Reference}}
				\\[1mm] \hline \hline
				Eg1 & Volterra 
				& $\mmatP{\phantom{-}ax(1-y)\\ -by(1-x)}$ 
				&  $\mmatP{a\\ b}=\mmatP{2\\ 1}$
				& $Box_{(1,3)}( 0.1)$
				& \cite{moore:diffEqn:09},
				\cite[p.13]{bunger:taylorODE:20}
				\\[2mm] \hline 
				Eg2 & Van der Pol 
				& $\mmatP{y \\ -c(1-x^2)y -x}$
				& $c=1$
				& $Box_{(-3,3)}( 0.1)$
				& \cite[p.2]{bunger:taylorODE:20}
				\\[1mm] \hline
				Eg3 & Asymptote
				& $\mmatP{x^2 \\ -y^2 + 	7x}$
				& N/A
				& $ Box_{(-1.5,8.5)}(0.01 )$
				& N/A
				\\[1mm] \hline
				Eg4 & Lorenz
				& $\mmatP{\sigma(y-x)\\ x(\rho-z)-y\\ xy-\beta z}$
				&
				$\mmatP{\sigma\\\rho\\\beta}=\mmatP{10\\ 28\\8/3}$
				& $ Box_{(15,15,36)}(0.001)$
				& \cite[p.11]{bunger:taylorODE:20}
				\\[1mm] \hline
				Eg5 & R\"ossler
				& $(-y-z,x+ay,b+z(x-c))$
				&
				$(a,b,c)=(0.2,0.2,5.7)$
				& $ Box_{(1,2,3)}(0.1)$
				& \cite{rossler:chaos:76}
				\\[1mm] \hline
		}}		
	\end{table}
	
\ssect{One Step Comparison of Taylor Tube formulas (\refTab{deltafull})}
	
	In \refTab{deltafull}, for each example and for
	any given admissible $(E_0, H, F_1)$
		(which is obtained from \stepA),
	we first compute
		$\olh = \min(H,h^{taylor}(H, \olM, \mu, \delta))$
	(see \refeq{h1}).
	Then we compute the full-enclosure of
		$\ivp(E_0,\olh)$ using, respectively, the formulas
	\refeq{direct} and \refeq{E1h1}.
	We compare the relative effectiveness of the 2 formulas
	using the ratio
	\beql{sigma}
		\sigma\as \frac{
			\displaystyle vol\!\left(
			\Direct(E_0,[0,\olh],F_1) \right)
		}{
			\displaystyle vol\!\left(
			E_1([0,\olh]) \right)
		},
	\eeql
	where $vol(\cdot)$ is the volume of a box.
	Thus $\sigma>1$ implies an improvement due to Taylor tube.
	In the table, we choose different values of
	$\delta$ and $k$ (with $p=k-1$).

	\begin{table*}
	{\tiny
	\centering
	\caption{
		Comparison of Full-Enclosures from $Direct(E_0,[0,h],F_1)$
		and $E_1([0,h])$.
	}
	\label{tab:deltafull}

	\begin{tabular}{c|>{\centering\arraybackslash}
		p{3cm}|c|>{\centering\arraybackslash}p{3.5cm}|c|c|c|c|c}
	\hline
	Eg & $E_0$ & $H$ & $F_1$  & $\delta$ & order
		& $\olh$ & $\mu$ & $\sigma$ \\
	\hline
	
	\multirow{12}{*}{Volterra}
	& \multirow{12}{3.8cm}{$(1.0,3.0)\pm(0.1,0.1)$}
	& \multirow{12}{*}{0.1}
	& \multirow{12}{4.2cm}{$(0.75,3.0)\pm(0.36,0.20)$}
	& \multirow{4}{*}{0.1}  & 1  & 0.0144  & 0.1378 & 1.00    \\
	& & & &                 & 3  & 0.1     & 0.4386 & 1.03 \\
	& & & &                 & 7  & 0.1     & 0.4171 & 1.07 \\
	& & & &                 & 20 & 0.1     & 0.4169 & 1.07 \\
	\cline{5-9}
	& & & & \multirow{4}{*}{0.01} & 1  & 0.00201 & 0.1048 & 1.09  \\
	& & & &                  & 3  & 0.0525  & 0.2398 & 1.30 \\
	& & & &                  & 7  & 0.1     & 0.36   & 1.69 \\
	& & & &                  & 20 & 0.1     & 0.36   & 1.69 \\
	\cline{5-9}
	& & & & \multirow{4}{*}{0.001} & 1  & 0.00020 & 0.1005 & 1.00    \\
	& & & &                   & 3  & 0.0275  & 0.17   & 1.23 \\
	& & & &                   & 7  & 0.1     & 1.9    & 1.85 \\
	& & & &                   & 20 & 0.1     & 0.35   & 1.85 \\
	\hline
	
	\multirow{12}{*}{Van der Pol}
	& \multirow{12}{3.8cm}{$(-3.0,3.0)\pm(0.1,0.1)$}
	& \multirow{12}{*}{0.05}
	& \multirow{12}{4.2cm}{$(-2.9,2.4)\pm(0.19,0.71)$}
	& \multirow{4}{*}{0.1}  & 1  & 0.0035  & 0.6314 & 1.10   \\
	& & & &                 & 3  & 0.0481  & 1.1    & 1.69 \\
	& & & &                 & 7  & 0.05    & 1.2    & 1.76 \\
	& & & &                 & 20 & 0.05    & 1.2    & 1.76 \\
	\cline{5-9}
	& & & & \multirow{4}{*}{0.01} & 1  & 0.00041 & 0.60   & 1.22   \\
	& & & &                  & 3  & 0.0239  & 0.83   & 1.34 \\
	& & & &                  & 7  & 0.05    & 1.1    & 2.02 \\
	& & & &                  & 20 & 0.05    & 1.2    & 2.02 \\
	\cline{5-9}
	& & & & \multirow{4}{*}{0.001} & 1 
			& $4.1\times 10^{-5}$ & 0.60   & 1.05    \\
	& & & &                   & 3  & 0.0114  & 0.70   & 1.16 \\
	& & & &                   & 7  & 0.05    & 1.0    & 2.09 \\
	& & & &                   & 20 & 0.05    & 1.1    & 2.09 \\
	\hline
	
	\multirow{12}{*}{Asymptote}
	& \multirow{12}{2.8cm}{\centering \pmcell{(-1.50,8.50)}{(0.01,0.01)}}
	& \multirow{12}{*}{0.04}
	& \multirow{12}{4.2cm}{$(-1.5,6.8)\pm(0.059,1.7)$}
	& \multirow{4}{*}{0.1}  & 1  & 0.00117 & 0.1167 & 1.02    \\
	& & & &                 & 3  & 0.0246  & 2.3    & 1.42 \\
	& & & &                 & 7  & 0.04    & 3.5    & 1.97 \\
	& & & &                 & 20 & 0.04    & 3.5    & 1.97 \\
	\cline{5-9}
	& & & & \multirow{4}{*}{0.01} & 1  & 0.00012 & 0.030  & 1.01    \\
	& & & &                  & 3  & 0.0116  & 1.0    & 1.22 \\
	& & & &                  & 7  & 0.04    & 3.4    & 2.14 \\
	& & & &                  & 20 & 0.04    & 3.4    & 2.15 \\
	\cline{5-9}
	& & & & \multirow{4}{*}{0.001} & 1 
			& $1.2\times 10^{-5}$ & 0.021  & 1.00    \\
	& & & &                   & 3  & 0.0054  & 0.47   & 1.10 \\
	& & & &                   & 7  & 0.0305  & 2.8    & 1.82 \\
	& & & &                   & 20 & 0.04    & 3.3    & 2.24 \\
	\hline
	
	\multirow{12}{*}{Lorenz}
	& \multirow{12}{2.8cm}
		{\centering \pmcell{(15.000,15.000,36.000)}{(0.001,0.001,0.001)}}
	& \multirow{12}{*}{0.027}
	& \multirow{12}{2.2cm}{\centering \pmcell{(15,13,38)}{(0.26,2.2,1.8)}}
	& \multirow{4}{*}{0.1}  & 1  & 0.000508 & 2.1   & 1.000 \\
	& & & &                 & 3  & 0.0151   & 2.9   & 1.000 \\
	& & & &                 & 7  & 0.027    & 3.3   & 1.000 \\
	& & & &                 & 20 & 0.027    & 3.3   & 1.000 \\
	\cline{5-9}
	& & & & \multirow{4}{*}{0.01} & 1 
			& $5.3\times 10^{-5}$ & 2.1   & 1.000 \\
	& & & &                  & 3  & 0.00719  & 2.5   & 1.000 \\
	& & & &                  & 7  & 0.027    & 3.3   & 1.000 \\
	& & & &                  & 20 & 0.027    & 3.3   & 1.000 \\
	\cline{5-9}
	& & & & \multirow{4}{*}{0.001} & 1 
			& $5.4\times 10^{-6}$ & 2.1   & 1.000 \\
	& & & &                   & 3  & 0.00387  & 2.3   & 1.000 \\
	& & & &                   & 7  & 0.0248   & 3.2   & 1.001 \\
	& & & &                   & 20 & 0.027    & 3.3   & 1.001 \\
	\hline
	
	\multirow{12}{*}{ R\"ossler}
	& \multirow{12}{2.8cm}
			{\centering \pmcell{(1.0,2.0,3.0)}{(0.1,0.1,0.1)}}
	& \multirow{12}{*}{0.1}
	& \multirow{12}{2.2cm}{\centering \pmcell{(0.74,2.1,2.3)}{(0.37,0.18,0.85)}}
	& \multirow{4}{*}{0.1}  & 1  & 0.00588 & 0.3155 & 1.00   \\
	& & & &                 & 3  & 0.1     & 0.71   & 1.62 \\
	& & & &                 & 7  & 0.1     & 0.70   & 1.73 \\
	& & & &                 & 20 & 0.1     & 0.70   & 1.73 \\
	\cline{5-9}
	& & & & \multirow{4}{*}{0.01} & 1  & 0.00063 & 0.30   & 1.00    \\
	& & & &                  & 3  & 0.0492  & 0.48   & 1.32 \\
	& & & &                  & 7  & 0.1     & 0.65   & 2.02 \\
	& & & &                  & 20 & 0.1     & 0.65   & 2.02 \\
	\cline{5-9}
	& & & & \multirow{4}{*}{0.001} & 1  & $6.3\times 10^{-5}$ & 0.29   & 1.00   \\
	& & & &                   & 3  & 0.0257  & 0.39   & 1.17 \\
	& & & &                   & 7  & 0.1     & 0.64   & 2.12 \\
	& & & &                   & 20 & 0.1     & 0.64   & 2.12 \\
	\hline\\
	
	\end{tabular}}
	\end{table*}
	
	From Table~\ref{tab:deltafull}, we see that $\sigma$ (last column)
	is generally $>1$ implying the effectiveness of Taylor tube.
	Moreover, increasing the Taylor
	degree causes the step size $h\taylor_p$ to increase
	(see the $\olh$ column).
	However, when $\olh$ is capped at the value $\olh=H$,
	higher Taylor degrees do lead to larger
	values of $\sigma$ (though it may be invisible in
	the table where we only show 2 digits).
	\ignore{ 
	When $\delta$ decreases, the step size $h\taylor_p$
	also decreases; however, the obtained full enclosure becomes tighter.
	Finally, from the last column we observe that $\sigma$ is typically
	greater than $1$. This shows that the Taylor-tube method effectively
	reduces the size of the full enclosure.
	}%
	
\ssect{Global Experiment (\refTab{tube-degree-combined})}
	Table~\ref{tab:tube-degree-combined}
	shows the running times of our End-Cover Algorithm on each
	example using varying Taylor degrees
		$(1, 3, 5, 7, 19)$.  The order $k$
	is fixed at $20$.  We also record the size $|\calC|$
	of the end-cover.  But in each case, we show two numbers,
	representing the standard End-Cover Algorithm using \Bisect\
	(\withB) and its variant (\NoB) in which
	the \Bisect\ subroutine is turned off.
	This is because the
	algorithm invokes \Bisect\ and \taylortube\
	in tandem, their interaction is quite complex.
	This allows us to study the effects of Taylor degree in isolation.
	 
	For each example, we add a special row whose tube
	is degree ``$19$ (\NoTube)''.  This is our End-Cover Algorithm
	with \taylortube\ turned off.  
	Comparing the row for degree $19$ with this special row
	is extremely insightful: we get a triple of timing in seconds.
	For example, for Volterra, the triple is
		$(21.7s, \cored{7.06s}, 19.3s)$,
	corresponding to
		(\withT+\NoB, \cored{\withT+\withB}, {\NoTube+\withB}).
	Inevitably, the middle value is smallest, showing
	that the combination of \Bisect\ and \taylortube\ is
	better than either one alone.
	
	Observe that the use of
	\Bisect\ not only improves the overall efficiency, but also tends to
	reduce the size of the end-covers.
	(this is seen by comparing the numbers in the
	\NoB\ and \withB\ columns.)

	Remark: compared to other tube degrees, the timing for
	degree $19$ is artificially lower since it matches the
	order $k=20$ used in \stepB\ of~\refeq{E1h1}. This allows the
	Taylor-tube computation to reuse quantities 
	computed in \stepB\
	such as the point-error term. 

	\begin{table*}
	\centering
	\caption{Performance of End-Cover algorithm under different
			Taylor degrees.}
	\label{tab:tube-degree-combined}
	\scriptsize
	\begin{tabular}{c|c|c|c|c|c|c|c}
	\hline
	\multirow{2}{*}{\textbf{Model}} &
	\multirow{2}{*}{$T$} &
	\multirow{2}{*}{$\veps$} &
	\multirow{2}{*}{\textbf{Taylor degree}} &
	\multicolumn{2}{c|}{Cover size |$\calC$|} &
	\multicolumn{2}{c}{\textbf{Time (s)}} \\
	\cline{5-6}\cline{7-8}
	& & & & \textbf{NoB} & \textbf{B} & \textbf{NoB} & \textbf{B} \\
	\hline
	
	\multirow{6}{*}{Volterra}
	& \multirow{6}{*}{$2$} & \multirow{6}{*}{$0.01$}
	& 1  & 5968 & 1024 & 47.5  & 8.05 \\
	& & & 3  & 5758 & 1024 & 44.5  & 7.82 \\
	& & & 5  & 5758 & 1024 & 46.2  & 10.6 \\
	& & & 7  & 5719 & 1024 & 47.3  & 11.0 \\
	\cline{4-8}
	& & & 19 & 2461 & 1024 & 21.7  & \cored{7.06} \\
	& & & 19 (\NoTube) & - & 1024 & - & 19.3 \\
	\hline
	
	\multirow{6}{*}{Van der Pol}
	& \multirow{6}{*}{$2$} & \multirow{6}{*}{$0.01$}
	& 1  & 13219 & 4174 & ~75.42~ & 35.58 \\
	& & & 3  & 9196  & 4195 & 60.88 & 35.99 \\
	& & & 5  & 6385  & 4111 & 89.94 & 34.74 \\
	& & & 7  & 4789  & 4102 & 99.55 & 34.48 \\
	\cline{4-8}
	& & & 19 & 4732  & 4098 & 44.70 & \cored{34.10} \\
	& & & 19 (\NoTube) & - & 4096 & - & 35.28 \\
	\hline
	
	\multirow{6}{*}{Asymptote}
	& \multirow{6}{*}{$1$} & \multirow{6}{*}{$0.01$}
	& 1  & 4933 & 4096 & 53.05 & 52.08 \\
	& & & 3  & 4936 & 4096 & 56.65 & 53.64 \\
	& & & 5  & 4930 & 4096 & 60.55 & 53.90 \\
	& & & 7  & 4705 & 4096 & 62.13 & 53.65 \\
	\cline{4-8}
	& & & 19 & 4642 & 4096 & 62.61 & \cored{51.14} \\
	& & & 19 (\NoTube) & - & 4096 & - & 51.846 \\
	\hline
	
	\multirow{6}{*}{Lorenz}
	& \multirow{6}{*}{$1$} & \multirow{6}{*}{$1.0$}
	& 1  & 246 & 8 & 10.6 & 0.293 \\
	& & & 3  & 246 & 8 & 14.0 & 0.296 \\
	& & & 5  & 246 & 8 & 16.2 & 0.302 \\
	& & & 7  & 246 & 8 & 18.6 & 0.300 \\
	\cline{4-8}
	& & & 19 & 211 & 8 & 36.4 & \cored{0.274} \\
	& & & 19 (\NoTube) & - & 8 & - & 1.75 \\
	\hline
	
	\multirow{6}{*}{R\"ossler}
	& \multirow{6}{*}{$2$} & \multirow{6}{*}{$1.0$}
	& 1  & 8 & 8 & 0.146 & 0.140 \\
	& & & 3  & 8 & 8 & 0.165 & 0.166 \\
	& & & 5  & 8 & 8 & 0.188 & 0.185 \\
	& & & 7  & 8 & 8 & 0.199 & 0.210 \\
	\cline{4-8}
	& & & 19 & 8 & 8 & 0.138 & \cored{0.111} \\
	& & & 19 (\NoTube) & - & 8 & - & 0.370 \\
	\hline\\
	
	\end{tabular}
	\end{table*}

\subsection{Conclusion}
The general introduction of
Taylor Tube methods represents a new tool in our arsenal for validated
IVP computation. We have shown experimentally that
they are effective in our algorithm.


\appendix
\section{Appendix A: Proofs}
\paragraph{Proof of \refLem{taylorStep}.}
	\ \\
	For simplicity, assume that $h=m\olh$ for some integer $m$, and write
		\[ t_i=i\olh, \qquad i=0,\dots,m.  \]
	The proof below is easily adjusted if the last step size is $<\olh$.
	
	For each $i=0,\dots,m-1$, let
		\[ \bfq_i \as \Tay[p]{\bfq_0,\olh}(t_i), \qquad
			g_i \as \|\bfq_i-\bfx(t_i)\|  \]
	and
		\[ \bfx_i:[0,\olh]\to\RR^n \]
	denotes the exact solution of
		$ \bfx_i'(s)=\bff(\bfx_i(s))$, $\bfx_i(0)=\bfq_i$.
	By hypothesis, $\Image(\bfq_i,\olh)\subseteq B_1$.
	Since $B_1$ is a full enclosure of $(\bfq_0,h)$, we also have
		$\bfx([t_i,t_{i+1}])\subseteq B_1$.
	Hence, we can invoke \refLem{FundIneq} to
	bound the error between $\bfx$ and $\bfx_i$
	over their common time interval:
		\beql{error-i}
		\|\bfx_i(s)- \bfx(t_i+s)\|
				\le \|\bfq_i- \bfx(t_i)\| e^{\olmu s}
				= g_i e^{\olmu s}.
		\eeql
	Next, we bound the error
	$E(t)$ between the Taylor curve and $\bfx$ at 
	time $t\in [0,h]$ as follows: 
	{\small \beqarrys
	 E(t) &\as &\|\Tay[p]{\bfq_0,\olh}(t)-\bfx(t)\|  
			\\
		&=& \|\Tay[p]{\bfq_i}(s)-\bfx(t_i+s)\|
			\\ &&\text{\qquad(where $i=\floor{t/\olh}$,
					$t=s+t_i$ and $0\le s<\olh$)}\\
		&\le& \|\Tay[p]{\bfq_i}(s)-\bfx_i(s)\|
				+ \|\bfx_i(s)-\bfx(t_i+s)\|
			\\ &&\text{\qquad(by triangular inequality and
				$\bfq_i=\bfx_i(0)=\bfx(t_i)$)}\\
		&\le& \|\Tay[p]{\bfq_i}(s)-\bfx_i(s)\|
				+ g_i e^{\olmu s}
			\\ &&\text{\qquad(by \refeq{error-i})}\\
		&=& \|f\supn[p+1](\bfq_i) \xi^{p+1}\|
				+ g_i e^{\olmu s}
			\\ &&\text{\qquad(Taylor remainder, $0\le \xi\le s$)}\\
		&\le& \olM s^{p+1}
				+ g_i e^{\olmu s}
			\\ &&\text{\qquad($\bfq_i\in B_1$ and
					$\|\bff\supn[p+1](B_1)\|\le \olM$).}
	\eeqarrys}
	
	Let $\phi(s)$ be the last expression in this derivation:
		\beql{phii}
			\phi_i(s) \as \olM s^{p+1} + g_i e^{\olmu s}.\eeql
	Note that $\phi_i(s)$ is convex,
		being the sum of convex functions
			$\olM s^{p+1}$ and $g_i e^{\olmu s}$,
	and so its extremal values are obtained at $\phi_i(0)$
	and $\phi_i(\olh)$.  Let $\phi:[0,h]\to\RR$ where
		$\phi(t)\as \phi_i(s)$ if $t=s+t_i$, $i=\floor{s/\olh}$.
	Thus $\phi$ is an upper bound on the error function,
			$$E(t)\le \phi(t).$$
	By construction $\phi(t)$ is a continuous function
	whose extremal values are attained at $\phi(t_i)$
	($i=0\dd m$).  Next define
	\beql{Gi}
		G_i \as \clauses{0 & \rmif\ i=0,\\
			 			\olM \olh^{p+1} + G_{i-1} e^{\olmu \olh}
			 						&\rmif\ i\ge 1.}
	\eeql
	From \refeq{phii}, we see that
	$\phi(t_i)\le G_i$ for each $i$.
	We claim that
		\beql{Gi2}
			G_i = \clauses{\olM \olh^{p+1} &\rmif\ \olmu=0,\\ 
					\olM \olh^{p+1}\left(
						\frac{e^{i \olmu\olh}-1}{e^{\olmu\olh}-1}
						\right) &\rmif\ \olmu\ne 0.}
		\eeql
	The case $\olmu=0$ follows from the formula
	for $\olmu\ne 0$ by applying L'H\^{o}pital's rule.
	Observe that the sequence
		$G_i$ is non-decreasing with $i$ (even
		if $\olmu<0$).
	We prove \refeq{Gi2} by induction on $i$.
	The case $i=0$ is clear, and for $i\ge 1$: 
	\beqarrys
		G_i &=& \olM \olh^{p+1} + G_{i-1} e^{\olmu \olh}
					& \text{(from \refeq{Gi})}\\
			&=& \olM \olh^{p+1} +
				\olM \olh^{p+1}\left(
				\frac{e^{(i-1) \olmu\olh}-1}{e^{\olmu\olh}-1} \right)
							e^{\olmu\olh}
					& \text{(inductive hypothesis for $G_{i-1}$)}\\
		\ignore{
			&=& \olM \olh^{p+1} + \olM \olh^{p+1}\left(
				\frac{e^{i\cdot \olmu\olh}- e^{\olmu\olh}}
							{e^{\olmu\olh}-1} \right) \\
			&=& \olM \olh^{p+1}\left( 1+
				\frac{e^{i\cdot \olmu\olh}- e^{\olmu\olh}}
							{e^{\olmu\olh}-1} \right)\\
		}
			&\vdots& \\
			&=& \olM \olh^{p+1}\left(
				\frac{e^{i\cdot \olmu\olh}- 1}
					{e^{\olmu\olh}-1} \right)
					& \text{(as claimed).}
	\eeqarrys

	Finally, for any $\delta>0$,
	we find an upper bound on $\olh$ in order to
	ensure that the Taylor curve $\Tay[p]{\bfq_0,\olh}$
	stays inside the $\delta$-tube of $\bfx$.
	There are 3 cases, depending on $\olmu$:
	
	\dt{CASE $1$ ($\olmu> 0$)}:
		The sequence $G_i$ is increasing for $i=0\dd m$, and
		so it suffices to choose $\olh$ so that $G_m\le \delta$:
		\beqarrys
		 G_m 
			&=& \olM\olh^{p+1} \left( \frac{e^{m \olmu \olh}- 1}
					{e^{\olmu\olh}-1} \right) \\
			&=& \olM\olh^{p+1} \left( \frac{e^{\olmu h}- 1}
					{e^{\olmu\olh}-1} \right) 
						&\text{(since $h=\olh m$)}\\
			&\le& \olM\olh^{p+1} \left( \frac{e^{\olmu h}- 1}
					{\olmu\olh} \right)
				& \text{(since $e^x-1\ge x$ for $x\ge 0$)}\\
			&=& 
\olM\olh^{p}
				\left( \frac{e^{\olmu h}- 1}
					{\olmu} \right).
        \eeqarrys
To ensure $G_m\le \delta$,
it suffices to require
\beqarrys
\olM\olh^{p}
				\left( \frac{e^{\olmu h}- 1}
					{\olmu} \right)&\le&
                    \delta
                    \\
		\olh^p &\le& \frac{ \delta}{\olM}
					\frac{\olmu }{(e^{\olmu h} -1)}\\
		\olh &\le& \big(\frac{ \delta}{\olM}
					\frac{\olmu }{(e^{\olmu h} -1)}\big)^{1/p}.
		\eeqarrys

	\dt{CASE $0$ ($\olmu= 0$):}
		This case can be derived from CASE 1 by
		taking the limit as $\olmu\to 0$, and applying
		L'H\^{o}pital's rule.

	\dt{CASE $-1$ ($\olmu< 0$)}:
		In this case, the above derivation
		must use this bound $e^x-1< x+ \half x^2$ when $x<0$.
		We need some preliminary facts:
		\benum[(a)]
		\item Suppose $\delta>0$ and $B<B^+<0$.
			Then $\delta>\frac{A}{B^+}$ implies $\delta>\frac{A}{B}$.
		\NOignore{%
			If $A\ge 0$, the conclusion clearly holds.
			If $A<0$ then
				$$\delta>\frac{A}{B^+}=\big| \frac{A}{B^+}\big|
						> \big|\frac{A}{B}\big| =\frac{A}{B}$$
						since $|B^+|>|B|$.
			}
		\item Suppose $\olmu<0$ and $\olh<2/|\olmu|$,
			then $\olmu+\half \olmu^2 \olh <0$.
		\NOignore{%
				$\olmu+\half \olmu^2 \olh
					=\olmu( 1+ \half \olmu\olh) < \olmu$
					since $0> \half \olmu \olh >-1$.
			}
		\item If follows from (a) and (b)
		that to prove that $\delta \ge G_m$, it is
		sufficient to prove
			$$\delta \ge
					\olM \olh^{p+1}\left(
				\frac{e^{\olmu h}- 1}
					{\olmu\olh+\half(\olmu\olh)^2} \right)
				=
					\olM \olh^{p}\left(
				\frac{e^{ \olmu h}- 1}
					{\olmu+\half\olmu^2\olh} \right)
			$$

Rearranging the above inequality, we obtain
\[
\olh^p \le
\frac{\delta\bigl(\olmu+\half\olmu^2\olh\bigr)}
{\olM\,(e^{\olmu h}-1)}
=
\frac{\delta\olmu}{\olM\,(e^{\olmu h}-1)}
+
\frac{\delta\olmu^2}{2\olM\,(e^{\olmu h}-1)}\,\olh .
\]
Define
\[
A \as \frac{\delta\olmu}{\olM\,(e^{\olmu h}-1)},
\qquad
B \as -\,\frac{\delta\olmu^2}{2\olM\,(e^{\olmu h}-1)} .
\]
Since $\olmu<0$ and $e^{\olmu h}-1<0$, both $A$ and $B$ are positive.
Hence it is enough to prove
\[
\olh^p + B\olh \le A .
\]
Now the function $\olh^p+B\olh$ is increasing for $\olh\ge 0$ since $B$ is positive.
Therefore, if we choose
\[
\olh \le \min\!\left\{ \left(\frac{A}{2}\right)^{1/p},\, \frac{A}{2B} \right\},
\]
then
\[
\olh^p \le \frac{A}{2},
\qquad
B\olh \le \frac{A}{2},
\]
and thus
\[
\olh^p+B\olh \le A.
\]
This proves $\delta\ge G_m$.

Finally, note that
\[
\frac{A}{2B}
=
\frac{\delta\olmu/(\olM(e^{\olmu h}-1))}
{-\,\delta\olmu^2/(\olM(e^{\olmu h}-1))}
=
-\frac{1}{\olmu},
\]
so our choice also implies $\olh\le -1/\olmu$, hence
$\olmu+\half\olmu^2\olh<0$, as required.

		\eenum

\ignore{
\smallskip
\noindent
{\bf Case 1: $\olmu>0$.}
\smallskip
\noindent
{\bf Case 2: $\olmu=0$.}
From
\[
g_i\le h\,\olM\olh^p,
\]
it follows that $g_i\le\delta$ for all $i$ provided
\[
\olh\le \left(\frac{\delta}{\olM h}\right)^{1/p}.
\]
This is exactly the limit of the formula in \refeq{h1} as $\olmu\to 0^{+}$.

\smallskip
\noindent
{\bf Case 3: $\olmu<0$.}
Using
\[
e^{\olmu\olh}-1\le \olmu\olh+\frac12\olmu^2\olh^2,
\]
we obtain the sufficient condition
\[
\olh
\le
\left(
\frac{2\olmu\,\delta}
{\olM(e^{\olmu h}-1)-\olM^2\delta}
\right)^{1/p},
\]
which implies $g_i\le\delta$ for all $i$.

Therefore, in all cases,
\[
g_i\le\delta,
\qquad i=0,\dots,m.
\]

Finally, since $\phi_i$ is convex and both endpoint values satisfy
\[
\phi_i(0)=g_i\le \delta,
\qquad
\phi_i(\olh)\ge g_{i+1}\le \delta,
\]
we conclude that
\[
\|\Tay[p]{\bfq_0,\olh}(t_i+s)-\bfx(t_i+s)\|\le\delta
\qquad\text{for all } s\in[0,\olh].
\]
Since this holds for every $i=0,\dots,m-1$, the entire Taylor curve
$\Tay[p]{\bfq_0,\olh}(t)$ lies in the $\delta$-tube of $\bfx$.
}

\paragraph{Proof of \refThm{taylorTube}.}
\ \\
Let
\[
\bfx_0 = \ivp(\bfq_0,h),\qquad \bfx\in \ivp(E_0,h).
\]
We want to show that
\beql{tria}
\|\bfx(t)-\Tay[p]{\bfq_0,\olh}(t)\|_2
\le
\|\bfx(t)-\bfx_0(t)\|_2
+
\|\bfx_0(t)-\Tay[p]{\bfq_0,\olh}(t)\|_2
\le
r_0 e^{\olmu t}+\delta
\eeql
for all $t\in[0,h]$.

Hence all hypotheses of \refLem{taylorStep} are satisfied, and therefore
\beql{point1}
\|\Tay[p]{\bfq_0,\olh}(t)-\bfx_0(t)\|\le \delta
\qquad\text{for all } t\in[0,h].
\eeql

Because $\olmu$ is a lognorm-bound on $F_1$, \refLem{FundIneq} gives
\beql{range}
\|\bfx(t)-\bfx_0(t)\|_2
\le
e^{\olmu t}\,\|\bfp-\bfq_0\|_2
\le
r_0 e^{\olmu t}
\qquad\text{for all } t\in[0,h].
\eeql
Combining \refeq{point1} and \refeq{range}, we obtain \refeq{tria}.

\begin{enumerate}[(a)]
\item
Taking $t=h$ in \refeq{tria}, we get
\[
\|\bfx(h)-\Tay[p]{\bfq_0,\olh}(h)\|_2
\le
r_0 e^{\olmu h}+\delta.
\]
Hence
\[
\bfx(h)\in
\mathrm{Ball}_{\Tay[p]{\bfq_0,\olh}(h)}
\left(r_0 e^{\olmu h}+\delta\right).
\]
 This ball is an end-enclosure for
$\ivp(E_0,h)$.

\item
Let
\[
r \as \max_{h'\in[0,h]}\left(r_0 e^{\olmu h'}+\delta\right).
\]
Then by \refeq{tria} for every $t\in[0,h]$,
\[
\|\bfx(t)-\Tay[p]{\bfq_0,\olh}(t)\|_2
\le
r_0 e^{\olmu t}+\delta
\le r.
\]
Therefore,
\[
\bfx(t)\in
\mathrm{Ball}_{\Tay[p]{\bfq_0,\olh}(t)}(r)
\ib
\Tay[p]{\bfq_0,\olh}(t)\pm[-r,r]^n.
\]
Taking the union over $t\in[0,h]$, we obtain
\[
\bfx([0,h])
\ib
\Tay[p]{\bfq_0,\olh}([0,h])\pm[-r,r]^n.
\]
This proves that
\[
\Tay[p]{\bfq_0,\olh}([0,h])\pm[-r,r]^n
\]
is a full-enclosure for $\ivp(E_0,h)$.
\end{enumerate}
\qed

\icmsOrArxiv{
}{ 
	
}%

\renewcommand{\alt}[2]{#2} 
\renewcommand{\alt}[2]{#1} 

\bibliographystyle{splncs04}
\alt{

}{
	\bibliography{st,yap,exact,com,gen,geo,alge,
		math,algo,mesh,tnt,ran,bk/bk,bk/bknew}
}%

\end{document}